\documentclass[10pt,conference]{IEEEtran}
\usepackage[utf8]{inputenc}
\usepackage{amsmath}
\usepackage{amsfonts}
\usepackage{amssymb}
\usepackage{graphicx}
\usepackage{cleveref}
\usepackage{relsize}
\usepackage{bbm}
\usepackage{enumerate}

\newcommand{\ie}{i.e.\ }
\newcommand{\eg}{e.g.\ }
\newcommand{\wrt}{w.r.t.\ }
\crefname{section}{\S}{\S\S}
\crefname{subsection}{\S}{\S\S}
\crefname{subsubsection}{\S}{\S\S}
\crefname{equation}{Eq.}{}
\crefname{figure}{Fig.}{}

\newcommand{\Pro}[1]{\mathbb{P}\left[ #1 \right]}

\newcommand{\round}{\mathrm{Round}}

\newcommand{\ceil}[1]{\lceil #1 \rceil}
\newcommand{\floor}[1]{\lfloor #1 \rfloor}

\newcommand{\fintvl}[1][x]{\mathlarger{\lfloor}#1,#1\mathlarger{\rceil}}
\newcommand{\Err}{\mathcal{E}}
\newcommand{\F}{\mathbb{F}}
\newcommand{\R}{\mathbb{R}}
\newcommand{\N}{\mathbb{N}}
\newcommand{\mop}{~\mathtt{op_m}~}
\newcommand{\iop}{~\mathrm{op}~}
\newcommand{\Exp}[1]{\mathbb{E}\left[#1\right]}
\newcommand{\one}{\mathbbm{1}}
\newcommand{\absv}[1]{\vert #1\vert}
\newcommand{\dt}{\frac{\partial}{\partial t}}

\usepackage[backgroundcolor=Apricot, textsize=footnotesize]{todonotes}

\title{A Probabilistic Approach to \\ Floating-Point Arithmetic}

\author{
\IEEEauthorblockN{Fredrik Dahlqvist}
\IEEEauthorblockA{Department of Electrical and\\ Electronic Engineering\\Imperial College London\\ f.dahlqvist09@imperial.ac.uk}\and 
\IEEEauthorblockN{Rocco Salvia}
\IEEEauthorblockA{School of Computing\\ University of Utah\\rocco@cs.utah.edu}\and 
\IEEEauthorblockN{George A. Constantinides}
\IEEEauthorblockA{Department of Electrical and\\ Electronic Engineering\\Imperial College London\\ g.constantinides@imperial.ac.uk}
}

\begin{document}
\maketitle

\begin{abstract}
Finite-precision floating point arithmetic unavoidably introduces rounding errors which are traditionally bounded using a worst-case analysis. However, worst-case analysis might be overly conservative because worst-case errors can be extremely rare events in practice. Here we develop a probabilistic model of rounding errors with which it becomes possible to estimate the likelihood that the rounding error of an algorithm lies within a given interval. Given an input distribution, we show how to compute the distribution of rounding errors. We do this exactly for low precision arithmetic, for high precision arithmetic we derive a simple approximation. The model is then entirely compositional: given a numerical program written in a simple imperative programming language we can recursively compute the distribution of rounding errors at each step of the computation and propagate it through each program instruction. This is done by applying a formalism originally developed by Kozen to formalize the semantics of probabilistic programs. We then discuss an implementation of the model and use it to perform probabilistic range analyses on some benchmarks.
\end{abstract}

\section{Introduction}

IEEE arithmetic \cite{ieee754} is traditionally modelled mathematically as follows \cite{higham2002accuracy}: if $x,y$ are two normal floating-point numbers and $\iop\in\{+,-,\times,\div\}$ is an infinite-precision arithmetic operation, then the floating-point precision implementation $\mop$ of $\iop$ must satisfy:
\begin{align}
x\mop y=(x\iop y)(1+\delta), \qquad\absv{\delta}\leq u\label{eq:traditional}
\end{align}
where $u$ is the unit roundoff for the given precision. \Cref{eq:traditional} says that the machine implementation of an arithmetic operation can induce a relative error of size $\delta$ \emph{for some} $\delta\in\left[-u,u\right]$. The `\emph{for some}' is essential: this is a \emph{non-deterministic model}, we have no control whatsoever over which $\delta$ appears in \cref{eq:traditional}. This means that numerical analysis based on this model must consider \emph{all} possible values $\delta$, \ie numerical analysis based on \cref{eq:traditional} is fundamentally a \emph{worst-case analysis}. 

It also follows from the perspective of \cref{eq:traditional} that any program doing arithmetic is, under this model, a non-deterministic program. Moreover, since the output of such a program might very well turn out to be the input of another program doing arithmetic, one should also consider non-deterministic inputs. This is precisely what happens in practice with tools for numerical analysis like the recent \cite{magron2017certified}, Daisy \cite{darulova2018daisy} or FPTaylor \cite{solovyev2018rigorous} which require for each variable of the program a range of possible values in order to perform a worst-case analysis.

However, for a wide variety of programs  it makes sense to assume that the inputs are \emph{probabilistic} rather than non-deterministic; that is to say we have some statistical model of the inputs of the program. This situation is in fact incredibly common. The inputs of one numerical routine are frequently generated randomly by another numerical routine, for example in a gradient descent optimization, a Bayesian inference algorithm, or a stochastic ray tracing algorithm. Similarly, sensors on a cyber-physical system can feed analog signals which are very well modelled statistically, to a numerical program processing these signals. 

If the inputs of a program have a known distribution, then it becomes possible, at least in principle, to ask the question: \textit{How likely are the inputs generating the worst-case rounding errors obtained from the non-deterministic model of \cref{eq:traditional}?} Typically, these inputs will occur very infrequently, and in this respect the non-deterministic model can be overly pessimistic since worst-case behaviours might be such rare events that they are never encountered in practice. 

In this paper we will explore a quantitative model which formally looks very similar to \cref{eq:traditional}, namely
\begin{align}
x\mop y=(x\iop y)(1+\delta), \qquad\delta\sim dist \label{eq:probabilistic}
\end{align}
but now $\delta$ is \emph{sampled} from $dist$, a probability distribution whose support is $\left[-u,u\right]$. In other words we move from a non-deterministic model of rounding errors to a \emph{probabilistic} model of rounding errors. This model will allow us to formalise and answer questions like: \textit{What is the distribution of outputs when rounding errors are taken into account?} \textit{What is the average rounding error?} \textit{What is the worst-case error with $99.9\%$ probability?}

As was mentioned above, model \eqref{eq:traditional} amounts to saying that any numerical program is a non-deterministic program. Completely analogously, in the perspective of \cref{eq:probabilistic} every numerical program is a \emph{probabilistic program}, that is to say a program which admits sampling as a native instruction. The study of probabilistic programs goes back to Kozen \cite{K81c} which modelled simple \texttt{while} programs containing an \emph{explicit} sampling instruction \texttt{random()}. In our setting any numerical program becomes a probabilistic program via an \emph{implicit} sampling operation which takes place whenever an arithmetic operation is performed. This implicit sampling is the only difference with the standard setting of \cite{K81c}, and we will otherwise understand how programs process randomness by following the framework laid out in \cite{K81c}. The study of probabilistic programs has recently witnessed a resurgence of interest driven by new applications in machine learning and statistical analysis of large datasets.

%
\emph{Related works: }
The probabilistic model of \cref{eq:probabilistic} is not new, it can be traced back to von Neumann and Goldstine \cite{von1947numerical} and is very similar to the so-called Monte-Carlo arithmetic of \cite{parker1997monte}. Within the signal processing 
community, simple probabilistic models of roundoff error are commonplace. Constantinides {\em et al.}~\cite{constantinides2004synthesis} study the propagation of fixed-point roundoff error through linear time-invariant computation, resulting in propagation of the
first and second statistical moments of the distributions.
 More recently, the model of \cref{eq:probabilistic} has been investigated by Higham \cite{higham2019new}
and Ipsen \cite{ipsen2019probabilistic}. Interestingly, because \cite{higham2019new} and \cite{ipsen2019probabilistic} are interested in large-dimensional problems, neither work needs to explicitly specify the distribution $dist$ in \cref{eq:probabilistic}. Instead, \cite{higham2019new} requires that each sample from $dist$ be independent and that $\Exp{\delta}=0$, whilst \cite{ipsen2019probabilistic} just requires that $\Exp{\delta}=0$. By using concentration of measure inequalities the authors then obtain probabilistic bounds which are independent of any particular choice of distribution. 
These bounds however are only applicable to a small class of problem (inner product, matrix multiplication and matrix factorisation algorithms) with very large inputs. Here we will derive a principled distribution $dist$ for the relative error $\delta$ and build a tool implementing the probabilistic model \cref{eq:probabilistic} to any small to medium-sized programs in a systematic way.

In \cite{probdaisy} the authors propose a hybrid approach: first, they
discretize input distributions and represent them as a dictionary which map each sub-interval to the corresponding probability (focal elements). On each sub-interval they combine probabilistic affine arithmetic, to propagate the error terms through the AST of the program, together with worst-case static analysis to bound the imprecision error term separately for each focal element. Their methodology results in a sound probabilistic error analysis, since the (unknown) error distribution is always bounded between proper upper and lower bounds by the analysis, but they inherit limitations of worst-case analysis in case of overflow and division by zero.\\
Our approach does not rely on worst-case error estimation, so it is a pure probabilistic model for the floating-point error of an expression.

\section{Two probabilistic models of rounding errors}\label{sec:rounding}

In order to use the probabilistic model given by \cref{eq:probabilistic} we must specify the distribution $dist$ of the random variable $\delta$. In this section we will show how to derive the distribution of rounding errors from first principles. This will yield a distribution which is computable for low precisions (\eg half-precision and lower) but becomes prohibitively expensive computationally for single- and double-precision. We will then show that very often the rounding error distribution can be approximated remarkably well by a simple distribution which we shall call the \emph{typical distribution}. The quality of this approximation increases with the working precision and we thus derive both an exact, computable model of rounding errors for low-precisions, and a simple but good approximating model of rounding errors for high-precisions.

\subsection{The exact rounding error distribution}\label{subsec:error_dist}

Conceptually, the key to our approach is to model the rounding operation probabilistically, \ie as an operation which adds a probabilistic relative error via
\begin{align}
x \longrightarrow x(1+\delta)\qquad \delta\sim dist.\label{eq:rounding}
\end{align}
Since each IEEE arithmetic operation can be understood as implicitly performing a rounding operation on the corresponding infinite-precision operation, the probabilistic rounding above naturally yields \cref{eq:probabilistic}. The key is thus to find a good candidate for the distribution $dist$ governing probabilistic rounding.

As discussed in the introduction, we consider numerical programs as probabilistic programs. In particular, all inputs come with probability distributions, and we should consider the variable $x$ in \cref{eq:rounding} as a sample from a real random variable $X$ with known probability distribution $\mathbb{P}$. It is then completely natural to require that:
\[
\frac{X-\round(X)}{X}~\sim~dist
\]
\ie $dist$ describes the distribution of the actual, deterministic rounding error of samples drawn from $X$. We will now explicitly compute $dist$. First we introduce some convenient notation. We define
\begin{align*}
\ceil{x}&\stackrel{\triangle}{=}\sup\{z\in\R\mid \round(z)=\round(x)\}\\
\floor{x}&\stackrel{\triangle}{=}\inf\{z\in\R\mid \round(z)=\round(x)\}.
\end{align*}
Whether $\ceil{x}$ is the maximal real which rounds to the same value as $x$, or just the supremum of this set, will in general depend both on $x$ and on the rounding convention, and similarly for $\floor{x}$. We also define the sets 
\[
\fintvl\stackrel{\triangle}{=}\left\{z\in\R\mid \round(z)=\round(x)\right\}
\]  
In particular if $z$ is a floating-point representable number -- notation $z\in\F$ -- then $\fintvl[z]$ is the collection of all reals rounding to $z$.

We choose to express the distribution $dist$ of relative errors in multiples of the unit roundoff $u$. This choice is arbitrary, but it allows us to normalize the distribution since the absolute value of the relative error is strictly bounded by $u$. In other words, we express the relative error as a distribution on $[-1,1]$ rather than $[-u,u]$. In order to compute the density function of $dist$ we proceed in the standard way by first computing the cumulative distribution function $c(t)$ and then taking its derivative. We therefore start by computing
\begin{align*}
c(t)\stackrel{\triangle}{=}&~\Pro{\frac{X-\widehat{X}}{X}\leq tu}\\
=&~\Pro{~\bigvee_{z\in\F}\left(\frac{X-z}{X}\leq tu\wedge X\in \fintvl[z]\right)}
\end{align*}
We now need to consider three special cases:
\begin{enumerate}
\item If $X\in \fintvl[0]$ then $\frac{X-\widehat{X}}{X}=1$ and thus (since $tu<1$):
\begin{align}\label{eq:cdfnot0}
\Pro{\frac{X-0}{X}\leq tu\wedge X\in \fintvl[0]}=0
\end{align}
\item If $X\in \fintvl[-\infty]$ then $\frac{X-\widehat{X}}{X}=\infty$ and thus 
\begin{align}\label{eq:cdfnotminf}
\Pro{\frac{X+\infty}{X}\leq tu\wedge X\in \fintvl[-\infty]}=0
\end{align}
\item Finally, if $X\in \fintvl[\infty]$ then $\frac{X-\widehat{X}}{X}=-\infty$ and thus 
\begin{align}\label{eq:cdfpinf}
\hspace{-3em}\Pro{\frac{X-\infty}{X}\leq tu\wedge X\in \fintvl[\infty]}=\Pro{X\in \fintvl[\infty]}
\end{align}
\end{enumerate}
Using the fact that \eqref{eq:cdfnot0}-\eqref{eq:cdfpinf} yield expressions which are independent of $t$ we get the density
\begin{align*}
d(t)=&\dt c(t)\\
=&\dt\sum_{z\in\F\setminus\{-\infty,0,\infty\}}\hspace{-1em}\Pro{\frac{X-z}{X}\leq tu\wedge X\in \fintvl[z]}\\
=& \sum_{z\in\F^-\setminus\{-\infty,0\}}\dt\Pro{\frac{z}{1-tu}\geq X\wedge X\in \fintvl[z]}+\\
&\sum_{z\in\F^+\setminus\{0,\infty\}}\dt\Pro{\frac{z}{1-tu}\leq X \wedge X\in \fintvl[z]}
\end{align*}
where $\F^+$ and $\F^-$ denote the positive (resp. negative) floating-point representable numbers.
Suppose now that $X$ is described by a probability density function $f:\R\to\R$, we then get:
\begin{align}
d(t)
=&\hspace{-1ex}\sum_{z\in\F^-\setminus\{-\infty,0\}}\dt\one_{\fintvl[z]}\left(\frac{z}{1-tu}\right) \int^{\frac{z}{1-tu}}_{\floor{z}} f(s)~ds+\nonumber 
\\
&\sum_{z\in\F^+\setminus\{0,\infty\}}\dt\one_{\fintvl[z]}\left(\frac{z}{1-tu}\right) \int^{\ceil{z}}_{\frac{z}{1-tu}} f(s)~ds\nonumber 
\\
=&\hspace{-2ex}\sum_{z\in\F^-\setminus\{-\infty,0\}}\hspace{-1ex}\one_{\fintvl[z]}\hspace{-4pt}\left(\frac{z}{1-tu}\right) f\hspace{-3pt}\left(\frac{z}{1-tu}\right) \frac{-uz}{(1-tu)^2}\nonumber
\\
& +\hspace{-2ex}\sum_{z\in\F^+\setminus\{0.\infty\}}\hspace{-2ex}\one_{\fintvl[z]}\hspace{-4pt}\left(\frac{z}{1-tu}\right) f\hspace{-3pt}\left(\frac{z}{1-tu}\right) \frac{uz}{(1-tu)^2}\nonumber 
\\
=&\hspace{-2ex}\sum_{z\in\F\setminus\{-\infty,0,\infty\}}\hspace{-2ex}\one_{\fintvl[z]}\hspace{-4pt}\left(\frac{z}{1-tu}\right) f\hspace{-3pt}\left(\frac{z}{1-tu}\right) \frac{u\absv{z}}{(1-tu)^2}\label{eq:errorDensity}
\end{align}
where $\one_{A}(x)$ is the usual indicator function which returns 1 if $x\in A$ and 0 otherwise. For low precisions, that is to say up to half-precision (5 bits exponent, 10 bits mantissa), it is perfectly possible to explicitly go through all floating point numbers and compute the density of the rounding error distribution $dist$ by using \cref{eq:errorDensity}. However this rapidly becomes too computationally expensive for higher-precision (since the number of floating-point representable numbers grows exponentially).

\subsection{The typical rounding error distribution}
\begin{figure}[h!]
\hspace{-1ex}\includegraphics[scale=0.60]{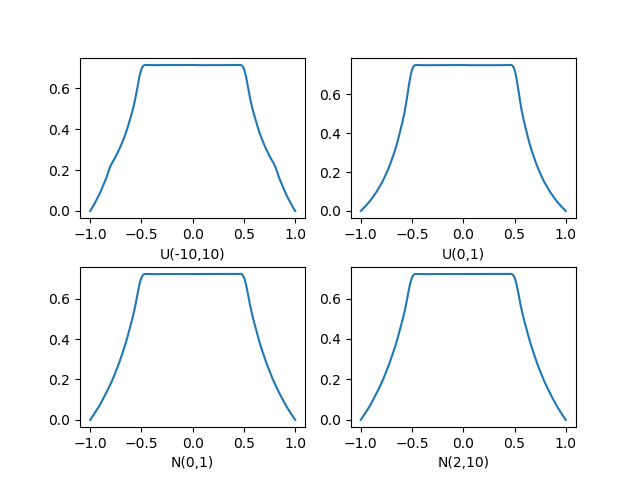}
\caption{Half-precision relative error distribution for four typical input distributions}
\label{fig:errdist}
\end{figure}

Interestingly, when computing the error density \cref{eq:errorDensity} for a wide variety of well-known input distribution, one very often obtains more or less the same curve. This phenomenon is illustrated in \cref{fig:errdist} where the half-precision error density computed via \cref{eq:errorDensity} is displayed for uniform distributions over $\left[-10,10\right]$ and $\left[0,1\right]$ and for normal distributions with parameters $\mu=0,\sigma=2$ and $\mu=2,\sigma=10$ respectively. The reader will notice immediately that all the curves are nearly identical.
In this section we will sketch how, under some regularity assumptions about the input distribution, the error density of \cref{eq:errorDensity} can be approximated by a simple, piecewise polynomial curve which we shall call the \emph{typical error distribution}. The precise mathematical derivation of this curve being relatively long, we refer the reader to Appendix for the full details. Here we will simply specify conditions under which the error distribution given by \cref{eq:errorDensity} is close to the typical distribution. The key observations which we shall make are that (a) the quality of the approximation increases with the working precision (henceforth $p$), and (b) the likelihood of the assumptions being satisfied also increases with $p$.

Let $z(e,s,k)$ denote the floating-point representable real with exponent $e$, sign $(-1)^s$ and mantissa $k$, and let $e_{min}$ and $e_{max}$ denote the smallest and largest exponents respectively. Given a $t\in \left[-1,1\right]$, one can show that the mantissas such that
\begin{align}
\one_{\fintvl[z(e,s,k)]}\left(\frac{z(e,s,k)}{1-tu}\right)=1\label{eq:membership}
\end{align}
are given by 
\begin{align}
k\leq 2^p\left(\frac{1}{\absv{t}}-1\right)-\frac{1}{2}\label{eq:mantissas}
\end{align}
Note that for $t\in \left[-\frac{1}{2},\frac{1}{2}\right]$ \cref{eq:mantissas} always holds, \ie all mantissa are compatible with \cref{eq:membership}.
We can now specify our assumptions.

\noindent \textbf{Assumption 0.} The probability density function is constant at the scale of the intervals between floating point numbers, \ie
\begin{align*}
\Pro{\round(x)=z}
&=\int_{\floor{z}}^{\ceil{z}} f(x)~dx 
\\
&\approx f\left(\frac{z}{1-tu}\right)(\ceil{z}-\floor{z})
\end{align*}
for all values of $t$ such that \cref{eq:membership} holds.

\noindent  \textbf{Assumption 1.} Writing $z(s,k)$ for $z(e_{min},s,k)$ assume that:
\[
\sum_{\substack{0\leq k<2^p\\ s\in\{0,1\}}}\hspace{-4pt}\one_{\fintvl[z]}\hspace{-4pt}\left(\frac{z(s,k)}{1-tu}\right)f\hspace{-3pt}\left(\frac{z(s,k)}{1-tu}\right)\hspace{-2pt}(\ceil{z(s,k)}-\floor{z(s,k)})\approx 0
\]
and similarly for $z(s,k)=z(e_{max},s,k)$.
Under Assumption 0, this condition says that the probability under $f$ of sampling a number whose rounding has exponent $e_{min}$ or $e_{max}$ is close to zero. This condition is certainly met for the distributions of \cref{fig:errdist} and half-precision.

\noindent \textbf{Assumption 2.} Given $t\in\left[-1,1\right]$, for every $k$ satisfying \cref{eq:mantissas} we assume
\[
\sum_{\substack{e_{min}<e<e_{max}\\ s\in\{0,1\}}}f\hspace{-2pt}\left(\frac{z(e,s,k)}{1-tu}\right)(\ceil{z(e,s,k)}-\floor{z(e,s,k)})\approx \frac{1}{2^p}
\]
Under assumption 0, this condition means that when rounding a sample drawn from the distribution $f$, all mantissas are equally likely.

\begin{figure}[ht!]
\includegraphics[scale=0.55]{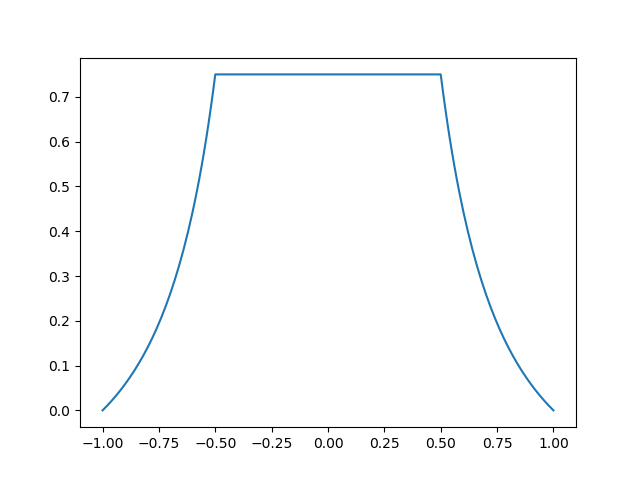}
\caption{Typical distribution of rounding errors (in unit roundoffs)}
\label{fig:typical}
\end{figure}

Under assumptions 0-2 one can show that for $t\in\left[-\frac{1}{2},\frac{1}{2}\right]$
\begin{align}
d(t)\approx\frac{1}{2^p(1-tu)^2}\left(\frac{2}{3}+\frac{3(2^{p}-1)}{4}\right)\label{eq:pdfmiddle}
\end{align}
Similarly, under assumptions 0-2 one has for $\absv{t}>\frac{1}{2}$:
\begin{align}
d(t)\approx\frac{1}{2^p(1-tu)^2}& \left(\frac{2}{3}+\frac{1}{2}\floor{2^p(\frac{1}{t}-1)-\frac{1}{2}}+\right. \nonumber 
\\
&\left.\frac{1}{2^{p+2}}(\floor{2^p(\frac{1}{t}-1)+\frac{1}{2}})^2\right)\label{eq:pdfwings}
\end{align}
where $\floor{2^p(\frac{1}{t}-1)+\frac{1}{2}}$ here denotes the usual floor function.
Combining \cref{eq:pdfmiddle} and \cref{eq:pdfwings} we get under assumptions 0-2 that as $p\to\infty$ the error density $d(t)$ is well approximated by the \emph{typical density}:
\begin{align}
d_{typ}(t)=\begin{cases}
\frac{3}{4}&\absv{t}\leq\frac{1}{2}
\\
\frac{1}{2}\left(\frac{1}{t}-1\right)+\frac{1}{4}\left(\frac{1}{t}-1\right)^2 & \absv{t}>\frac{1}{2}
\end{cases}\label{eq:typicalpdf}
\end{align} 
which is represented in \cref{fig:typical} and is clearly a good approximation of the exact densities of \cref{fig:errdist}.

\textbf{Remark:}
For well-behaved density functions, assumptions 0 and 1 become increasingly likely to hold as the precision $p$ and the exponent range $e_{min}$ to $e_{max}$ increase.

\section{Probabilistic interpretation of simple expressions}
In this section we present a class of numerical programs for which it is possible to \emph{formally} compute an output distribution given an input distribution, whilst implementing the probabilistic model of IEEE arithmetic operations given by \cref{eq:probabilistic}. In \cref{sec:exp} we will present some initial steps towards \emph{numerically} implementing the model presented in this section.

\subsection{A simple syntax}\label{subsec:syntax}

Our class of numerical program given by the following simple grammar.
\begin{align*}
&\textbf{Terms: }  
\\
& \tt t::= r \mid x_i \mid t\mop t \qquad \mathtt{r}\in\F, \mathtt{i}\in\N, \mop\hspace{-1ex}\in\{+.-,\times,\div\}
\\
&\textbf{Tests: } 
\\
& \tt b::= t < r \mid t > r \mid t == r  \qquad \mathtt{r}\in\F
\\
&\textbf{Programs:}
\\ 
&\tt p::= skip \mid \mathtt{x_i:=t}\mid \mathtt{p~;~p} \mid if~b~then~p~else~p
\end{align*}

For every expression $\mathtt{p}$, we consider the list $(\mathtt{x_1,\ldots,x_n})$ of variables appearing in $\mathtt{p}$. We view all variables as public variables and we associate with the list $(\mathtt{x_1,\ldots,x_n})$ a multivariate random variable (random vector) $\mathbf{X}$ modelling its (probabilistic) state. We will call this data the \emph{probabilistic context} and denote it $\{(\mathtt{x_1},\ldots,\mathtt{x_n})\sim\mathbf{X}\}$. We will denote by $X_i$ the $i^{th}$ marginal of $\mathbf{X}$. The cumulative distribution function of $\mathbf{X}$ will be given by
\[
\Pro{X_1\leq x_1,\ldots, X_n\leq x_n}.
\] 
We assume a fixed exponent range $e_{min}, e_{max}$ and precision level $p$ throughout.

\subsection{Random variable arithmetic}\label{subsec:rvarithmetic}

Whilst the probabilistic interpretation of programs will be defined in terms of an operator updating the probabilistic context, \ie an operator sending random vectors to random vectors, the probabilistic interpretation of \emph{terms} will be defined in terms of arithmetic operation on (univariate) random variables. We briefly review arithmetic operations on random variables which posses a density function (\wrt the Lebesgue measure) translate into operations on these densities \cite{springer1979algebra}. In particular the density of the sum of two independent random variables is given by the convolution of the densities. In more detail one has the following correspondence:
\begin{align}
X+Y&\sim f_X\oplus f_Y(t)=\int_{-\infty}^{\infty} f_X(x)f_Y(t-x)~dx\label{eq:pdfplus}\\
X-Y&\sim f_X\ominus f_Y(t)=\int_{-\infty}^{\infty} f_X(x)f_Y(x-t)~dx\label{eq:pdfminus}\\
X\times Y&\sim f_X\otimes f_Y(t)=\int_{-\infty}^{\infty} \frac{1}{\absv{x}}f_X(x)f_Y\left(\frac{t}{x}\right)dx\label{eq:pdftimes}\\
X\div Y&\sim f_X\oslash f_Y(t)=\int_{-\infty}^{\infty} \absv{x}f_X(x)f_Y(tx)dx\label{eq:pdfdiv}
\end{align}
Similarly, addition and multiplication by a scalar correspond to
\begin{align*}
\alpha+X&\sim (\alpha\oplus f_X)(t)=f(t+\alpha)\\
\alpha X&\sim (\alpha\otimes f_X)(t)=\alpha f(\alpha t)
\end{align*}

\subsection{Probabilistic interpretation}\label{subsec:sem}

For any univariate random variable $X$ with density $f_X$, we define the \emph{error of $X$} -- notation $\Err(X)$ -- as the random variable $\frac{X-\widehat{X}}{X}$ whose density is given by \cref{eq:errorDensity}. For any $n$-dimensional multivariate random variable $\mathbf{X}$ we write $X_i, 1\leq i\leq n$ for its $i^{th}$ \emph{marginal} distribution. The probabilistic interpretation developed below dates back to  \cite{K81c}.

\subsubsection{Input quantization} Given a probabilistic context, the first question is to decide whether we need to model an initial quantization. This would correspond for example to modelling the quantization of an analog input (\eg a sensor) or of an input generated at a higher precision level (\eg a routine at half-precision level receiving input in double-precision). If we choose to model this step, then we need to add a probabilistic error term to each input, this is achieved by the probabilistic quantization sending $\mathbf{X}$ to the multivariate distribution with cumulative distribution function
\[
\Pro{X_1(1+\Err(X_1))\leq x_1\ldots\leq X_n(1+\Err(X_n))}
\] 
where the explicit computation of the random variables $X_1(1+\Err(X_i))$ can be performed using the densities derived in \cref{sec:rounding} and the operations on densities defined in \cref{subsec:rvarithmetic}.

\subsubsection{Terms} will be modelled as (univariate) random variables using the following inductive definition:
\begin{align}
model(\mathtt{r})&= r\text{, the constant r.v.} \nonumber\\
model(\mathtt{x_i})&= X_i \nonumber\\
model(\mathtt{t_1\mop t_2})&=\left(model(\mathtt{t_1})\iop model(\mathtt{t_2})\right)\cdot  \nonumber\\
& \hspace{12pt}\left(1+\Err(model(\mathtt{t_1})\iop model(\mathtt{t_2}))\right)\label{eq:termsem}
\end{align}

Note how we model the arithmetic operations in accordance with the fundamental model of \cref{eq:probabilistic}: we first compute the random variable $model(\mathtt{t_1})\iop model(\mathtt{t_2})$ (typically using the operations on densities defined in \cref{subsec:rvarithmetic}) which corresponds to the infinite-precision operation, and we then add a probabilistic error term whose distribution is computed from the distribution of $model(\mathtt{t_1})\iop model(\mathtt{t_2})$ itself using the densities derived in \cref{sec:rounding}.

\subsubsection{Tests} correspond to the obvious subset of $\R^n$ generated by the comparisons. Thus if $\mathtt{t}(x_1,\ldots,x_n)$ is a term in $n$ variables 
\[
model(\mathtt{t<r})=\{(x_1,\ldots,x_n)\in\R^n\mid \mathtt{t}(x_1,\ldots,x_n)<r\}
\]
and similarly for $\mathtt{t>r}$ and $\mathtt{t==r}$. 

\subsubsection{Expressions} will be modelled as operations sending multivariate random variables to multivariate random variables. In effect, updating the probabilistic context.

\begin{enumerate}[(i)]
\item $\mathtt{skip}$:
\[
model(\mathtt{skip})(\mathbf{X})=\mathbf{X}
\]
\item Assignments: $model(\mathtt{x_i:=t})(\mathbf{X})$ is the multivariate random variable whose cumulative distribution function is given by
\begin{align*}
\mathbb{P}[&X_1\leq x_1,\ldots,X_{i-1}\leq x_{i-1},model(\mathtt{t})\leq x_i, 
\\
& X_{i+1}\leq x_{i+1},\ldots, X_{n}\leq x_n]
\end{align*}
\item Sequential composition:
\[
model(\mathtt{p_1~;~p_2})(\mathbf{X})=model(\mathtt{p_2})(model(\mathtt{p_1})(\mathbf{X}))
\]
\end{enumerate}

For $\mathtt{if~then~else}$ statements we need to introduce the following notation: if $B\subseteq \R^n$ is a measurable subset of $\R^n$ and $\mathbf{X}$ is a multivariate random variable in $R^n$, then $\mathbf{X}_B$ will denote the multivariate random variable with cumulative distribution function
\[
\Pro{X_1\leq x_1,\ldots,X_n\leq x_n \wedge (X_1,\ldots,X_n)\in B}
\]
With this notation we can now define
\begin{enumerate}
\item[(iv)] Conditionals:
\begin{align*}
 model&(\mathtt{if~b~then~p_1~else~p_2})(\mathbf{X})=
\\
&model(\mathtt{p_1})\hspace{-3pt}\left(\mathbf{X}_{model(\mathtt{b})}\right)+model(\mathtt{p_2})\hspace{-3pt}\left(\mathbf{X}_{model(\mathtt{b})^c}\right)
\end{align*}
where $model(\mathtt{b})^c$ is the complement of $model(\mathtt{b})$
\end{enumerate}

The rules given above allows us (in principle) to compute the output distribution of any expression given by the grammar of \cref{subsec:syntax}. By construction this output will include the cumulative and combined effect of every probabilistic rounding occurring through arithmetic operations via the model \eqref{eq:probabilistic}.

\section{Experimental results}\label{sec:exp}

We have implemented the probabilistic interpretation of the class of terms (defined in \cref{subsec:syntax}) whose syntax has the structure of a \emph{tree}, that is to say terms where variables are not repeated (in which case the syntax would be a DAG). We can thus compute the probabilistic interpretation of $\tt (x+y)/(z\ast t)$ but not of $\tt (x+y)/(x\ast y)$. The reason for this current limitation is that two occurrences of the same variable must be interpreted as two perfectly correlated random variables, \ie a probability distribution whose support lies on the diagonal of $\R^2$. Such a distribution cannot have a two-dimensional density function, and this renders their representation highly-non trivial. For the time being we leave this problem to future work.

As shown in \cref{subsec:sem}, the interpretation of terms is computed by performing arithmetic operations on random variables. For random variables whose probability distribution is representable by density functions, \cref{subsec:rvarithmetic} shows how these arithmetic operations can be implemented in practice. A sophisticated version of \cref{eq:pdfplus}-\cref{eq:pdfdiv} is implemented in the Python library PaCAL
\cite{jaroszewicz2012arithmetic}\cite{korzen2014pacal} on which we have based the evaluation of the probabilistic interpretation of terms.

\subsection{Rounding error distribution}
Our implementation focuses on computations performed in low precision, that is to say half-precision or lower. We therefore compute the distribution of relative errors using the exact analytic formula \cref{eq:errorDensity}.

PaCAL concretely implements random variables as a density function represented by a piecewise Chebyshev interpolation polynomial. We therefore use the same representation and concretely represent the density \cref{eq:errorDensity} as a Chebyshev interpolating polynomial. The polynomial interpolation is performed using the library \texttt{pychebfun}, a Python implementation of the library \texttt{chebfun} \cite{battles2004extension} developed to perform accurate and fast computations on functions represented as Chebyshev polynomial interpolations.

\begin{figure}[h!]
\includegraphics[scale=0.22]{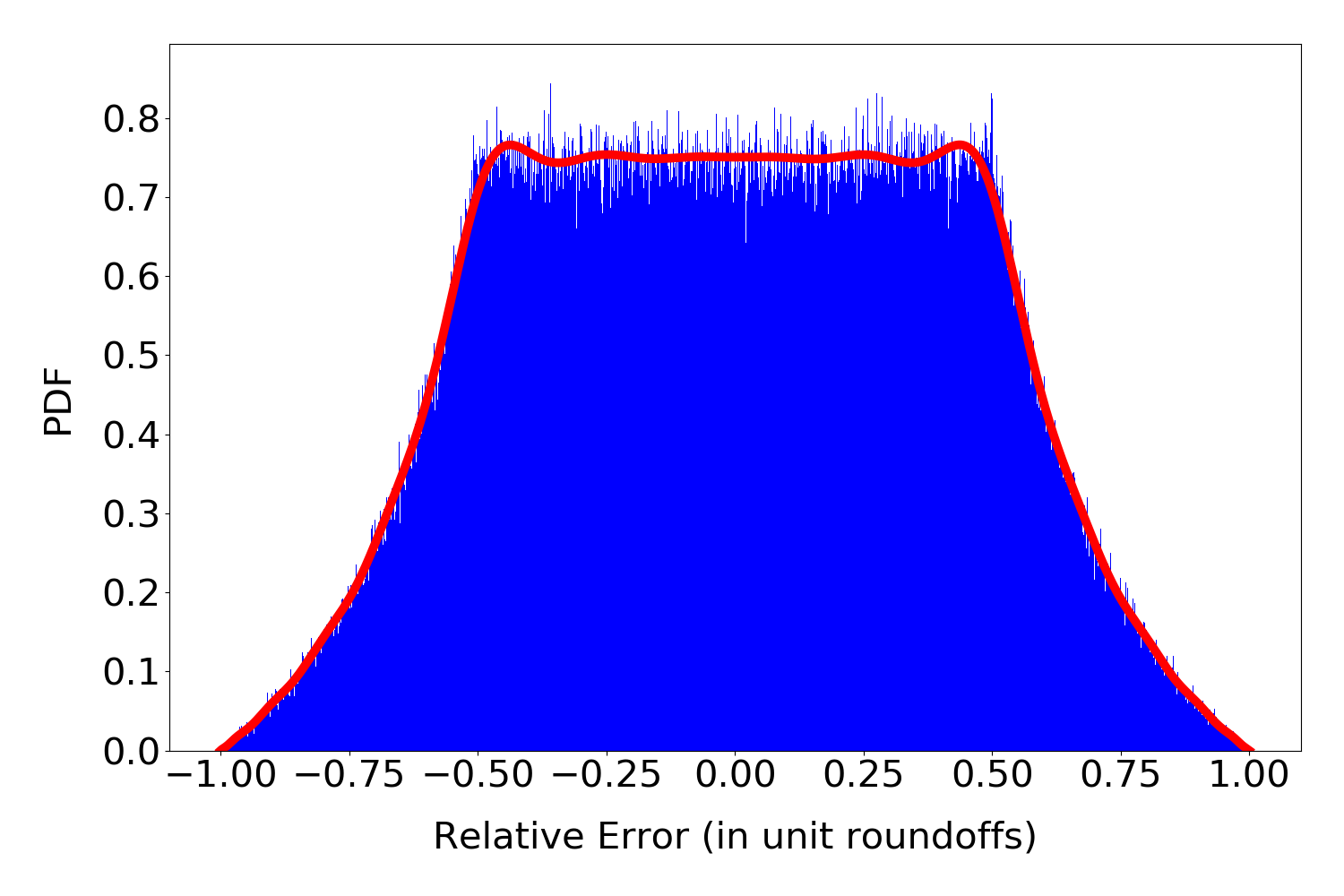}
\caption{Implemented density function vs Monte-Carlo simulation.}
\label{fig:impl:errdist}
\end{figure}

\cref{fig:impl:errdist} shows a plot of the half-precision error distribution for a random variable distributed uniformly on $\left[0,1\right]$. The red line is the density represented as a Chebyshev interpolating polynomial, the blue area is a histogram of the relative rounding error on one million samples. Note again the similarity with the typical density function represented in \cref{fig:typical}.

\subsection{Range analysis}
With a concrete representation of the error distribution associated with an input distribution we can recursively evaluate the probabilistic interpretation of terms following \eqref{eq:termsem}. This computes the final output distribution (including probabilistic rounding errors) of a term, given input distributions for each of its variables. A simple application of this output distribution is \emph{range analysis}. 

\begin{figure}[h!]
\includegraphics[scale=0.22]{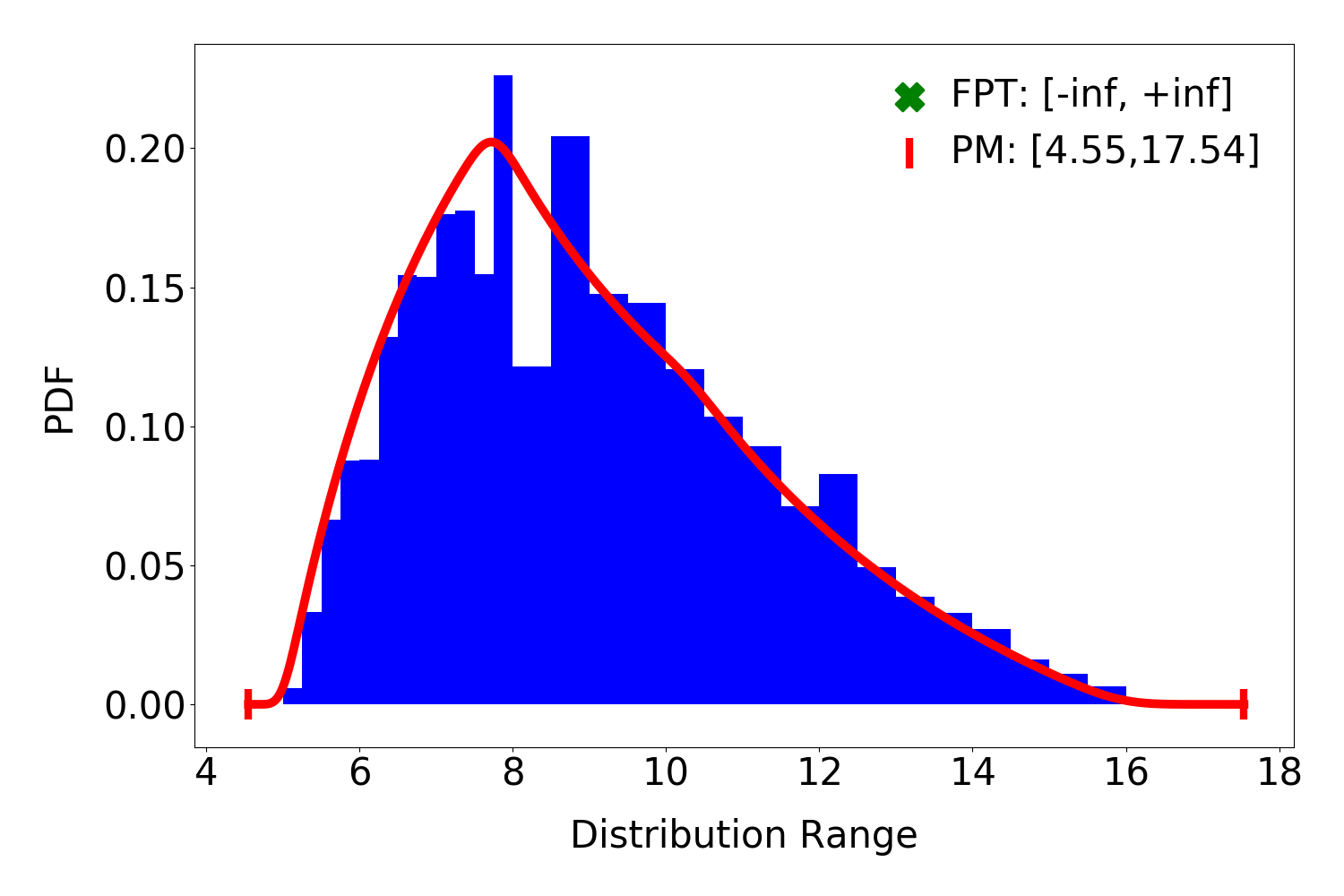}
\caption{Analytical and simulated output distribution for \texttt{x0/x1}.}
\label{fig:impl:div}
\end{figure}

We first consider a simple application of range analysis, namely the detection of overflows. To illustrate our probabilistic approach consider two variables $\tt x0, x1$ with $\tt x0$ distributed uniformly on $\left[10, 15.5\right]$ and $\tt x1$ distributed uniformly on $\left[0.97,2\right]$. Suppose that the working precision is 3 bits for the exponent and 3 bits for the mantissa and that we are interested in the term $\tt x0/x1$. In infinite precision $15.5/0.97 < 16$, the largest representable number at the given precision level, and there is no overflow. In reduced precision however overflow can occur. In this case, an analyser like FPTaylor will correctly detect the overflow and return an infinite range, but will not be able to quantify the \emph{likelihood of overflow}. \cref{fig:impl:div} shows the output distribution of the term $\tt x0/x1$ (red line), the support of the output distribution (`PM' in the legend), the output range of FPTaylor (`FPT' in the legend), and a histogram of the computation $\tt x0/x1$ performed on one million samples in reduced precision. Analytically, the probability of overflow is 0.0775\%, and empirically 0.0642\% of the samples overflow.

\begin{figure}[h!]
	\includegraphics[scale=0.22]{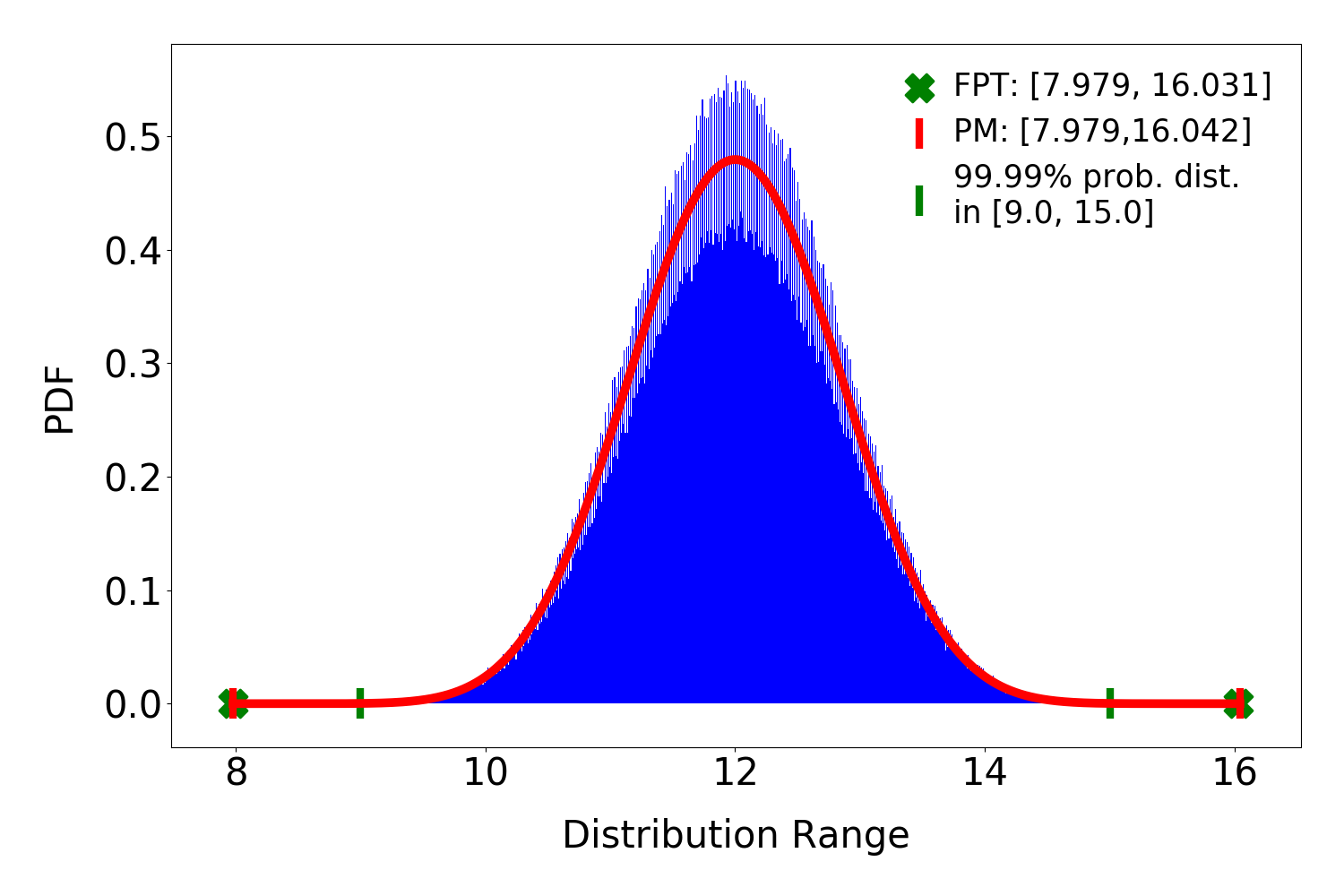}
	\caption{Analytical and simulated output distribution for \texttt{test02\_sum8less}.}
	\label{fig:impl:benchmark}
\end{figure}

Next, we consider two benchmarks from FPBench~\cite{fpbench} which our current setup can handle. In both cases we take \emph{half-precision} as the working precision. The first benchmark is \texttt{test02\_sum8less} which is given by the term: 
\[
\tt{((((((x0 + x1) + x2) + x3) + x4) + x5) + x6) + x7}
\]
with each $\tt xi$ assumed to be uniformly distributed on the interval $\left[1,2\right]$. The result of the range analysis are shown in \cref{fig:impl:benchmark}. The support of the output distribution provides marginally worse error bounds than FPTaylor and is in good agreement with one million Monte-Carlo simulations, that is to say one million evaluations of the term \texttt{test02\_sum8less} in reduced precision where is input value is sampled uniformly from $\left[1,2\right]$. Crucially, we can probabilistically tighten these bound: at 99.99\% confidence we can say that the output range lies in the interval $\left [9.0,15.0\right]$.

\begin{figure}[h!]
	\includegraphics[scale=0.177]{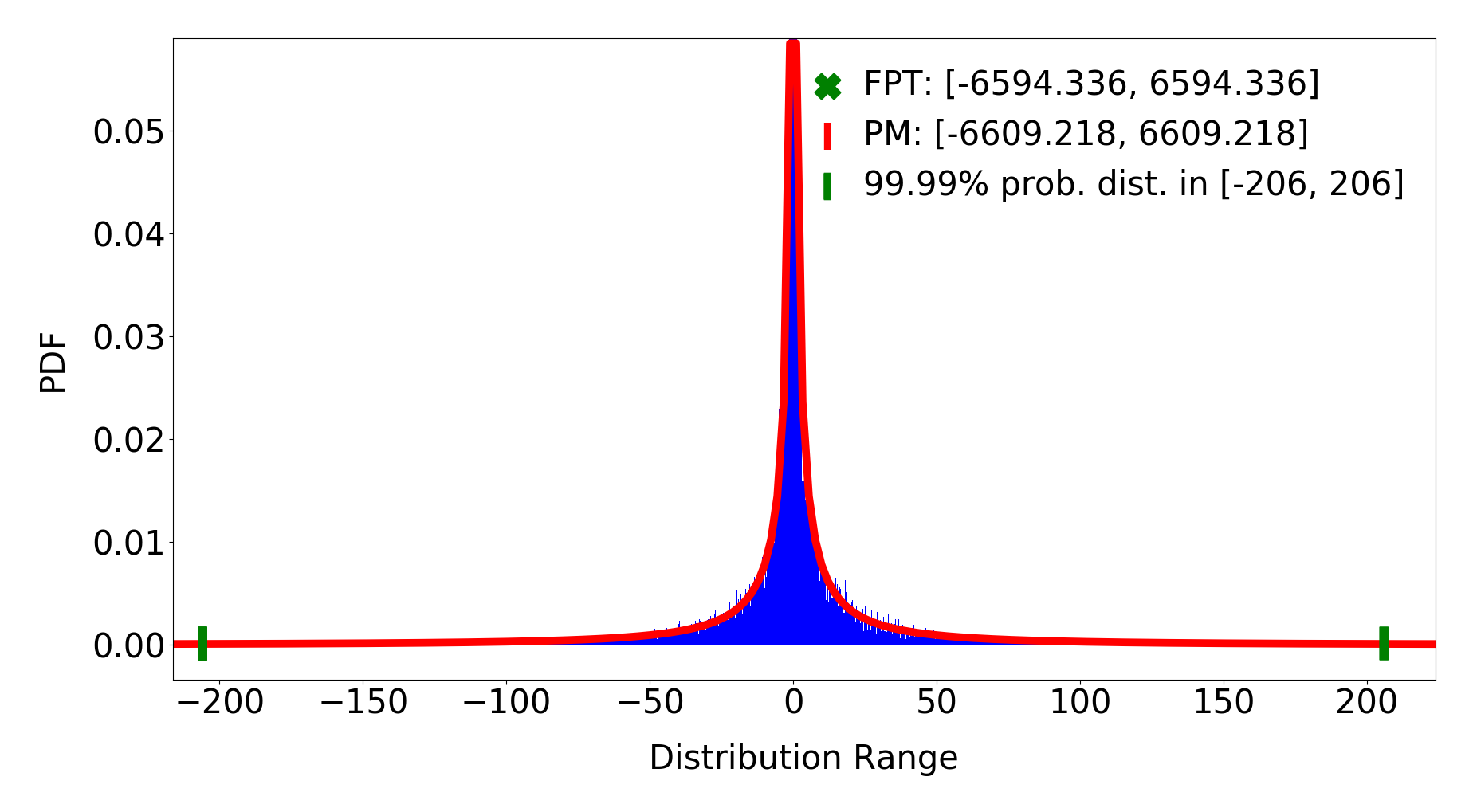}
	\caption{Analytical and simulated output distribution for \texttt{test02\_mul8less}}
	\label{fig:impl:benchmark2}
\end{figure}

The second benchmark provides an even stronger case for considering probabilistic range analysis. Consider the multiplicative cousin of \texttt{test02\_sum8less} given by the term \texttt{test02\_mul8less}:
\[
\tt{((((((x0 * x1) * x2) * x3) * x4) * x5) * x6) * x7}
\]
with each $\tt xi$ assumed to be uniformly distributed on the interval $\left[-3,3\right]$. The range analysis \texttt{test02\_mul8less} is displayed in \cref{fig:impl:benchmark2}. Again the support of the output distribution is marginally wider that the one provided by FPTaylor, but with 99.99\% confidence we can tighten the output range by about five orders of magnitude in base 2 to $\left[-206,206\right]$.

\bibliographystyle{plain}
\bibliography{constantinides-dahlqvist}

\begin{thebibliography}{10}

\bibitem{battles2004extension}
Zachary Battles and Lloyd~N Trefethen.
\newblock An extension of matlab to continuous functions and operators.
\newblock {\em SIAM Journal on Scientific Computing}, 25(5):1743--1770, 2004.

\bibitem{constantinides2004synthesis}
George~A Constantinides, Peter~YK Cheung, and Wayne Luk.
\newblock {\em Synthesis and optimization of DSP algorithms}.
\newblock Springer Verlag, 2004.

\bibitem{fpbench}
Nasrine Damouche, Matthieu Martel, Pavel Panchekha, Chen Qiu, Alexander
  Sanchez-Stern, and Zachary Tatlock.
\newblock Toward a standard benchmark format and suite for floating-point
  analysis.
\newblock In {\em International Workshop on Numerical Software Verification},
  pages 63--77. Springer, 2016.

\bibitem{darulova2018daisy}
Eva Darulova, Anastasiia Izycheva, Fariha Nasir, Fabian Ritter, Heiko Becker,
  and Robert Bastian.
\newblock Daisy-framework for analysis and optimization of numerical programs
  (tool paper).
\newblock In {\em International Conference on Tools and Algorithms for the
  Construction and Analysis of Systems}, pages 270--287. Springer, 2018.

\bibitem{higham2002accuracy}
Nicholas~J Higham.
\newblock {\em Accuracy and stability of numerical algorithms}, volume~80.
\newblock Siam, 2002.

\bibitem{higham2019new}
Nicholas~J Higham and Theo Mary.
\newblock A new approach to probabilistic rounding error analysis.
\newblock {\em SIAM Journal on Scientific Computing}, 41(5):A2815--A2835, 2019.

\bibitem{ipsen2019probabilistic}
Ilse~CF Ipsen and Hua Zhou.
\newblock Probabilistic error analysis for inner products.
\newblock {\em arXiv preprint arXiv:1906.10465}, 2019.

\bibitem{jaroszewicz2012arithmetic}
Szymon Jaroszewicz and Marcin Korze{\'n}.
\newblock Arithmetic operations on independent random variables: A numerical
  approach.
\newblock {\em SIAM Journal on Scientific Computing}, 34(3):A1241--A1265, 2012.

\bibitem{korzen2014pacal}
Marcin Korzen and Szymon Jaroszewicz.
\newblock {PaCAL}: A python package for arithmetic computations with random
  variables.
\newblock {\em Journal of Statistical Software}, 57(10):5, 2014.

\bibitem{K81c}
Dexter Kozen.
\newblock Semantics of probabilistic programs.
\newblock {\em J. Comput. Syst. Sci.}, 22(3):328--350, June 1981.

\bibitem{probdaisy}
Debasmita Lohar, Milos Prokop, and Eva Darulova.
\newblock Sound probabilistic numerical error analysis.
\newblock In {\em International Conference on Integrated Formal Methods}, pages
  322--340. Springer, 2019.

\bibitem{magron2017certified}
Victor Magron, George Constantinides, and Alastair Donaldson.
\newblock Certified roundoff error bounds using semidefinite programming.
\newblock {\em ACM Transactions on Mathematical Software (TOMS)}, 43(4):34,
  2017.

\bibitem{ieee754}
{Microprocessor Standards Committee of the IEEE Computer Society}.
\newblock {IEEE Standard for Floating-Point Arithmetic}, June 2019.

\bibitem{parker1997monte}
Douglass~Stott Parker.
\newblock {\em Monte Carlo Arithmetic: exploiting randomness in floating-point
  arithmetic}.
\newblock University of California (Los Angeles). Computer Science Department,
  1997.

\bibitem{solovyev2018rigorous}
Alexey Solovyev, Marek~S Baranowski, Ian Briggs, Charles Jacobsen, Zvonimir
  Rakamari{\'c}, and Ganesh Gopalakrishnan.
\newblock Rigorous estimation of floating-point round-off errors with symbolic
  taylor expansions.
\newblock {\em ACM Transactions on Programming Languages and Systems (TOPLAS)},
  41(1):20, 2018.

\bibitem{springer1979algebra}
M.D. Springer.
\newblock {\em The algebra of random variables}.
\newblock Probability and Statistics Series. Wiley, 1979.

\bibitem{von1947numerical}
John Von~Neumann and Herman~H Goldstine.
\newblock Numerical inverting of matrices of high order.
\newblock {\em Bulletin of the American Mathematical Society},
  53(11):1021--1099, 1947.

\end{thebibliography}

\appendix

\section*{Derivation of the typical distribution}

We start by examining the quantity
\[
\one_{\fintvl[z]}\left(\frac{z}{1-tu}\right)
\]
which appears in \cref{eq:errorDensity}. Given a $t\in\left[-1,1\right]$, we determine which floating point numbers $z\in\F\setminus\{-\infty,0,\infty\}$ make $\one_{\fintvl[z]}\left(\frac{z}{1-tu}\right)=1$. We will write floating point numbers $z\in\F$ in the form
\[
z=z(s,e,k):=(-1)^s ~ 2^e \left(1+\frac{k}{2^p}\right)
\]
where $s\in\{0,1\}$, $e$ ranges from $e_{min}=1-n$ to $e_{max}=n$ and $k$ ranges from 0 to $2^p-1$ (and represents the mantissa). We can explicitly compute $\floor{z}$ and $\ceil{z}$. When $s=0$ we have:
\begin{align*}
\floor{z}=
\begin{cases}
2^{e-1} & \text{if }s=0, e=e_{min}, k=0\\
2^{e-1}\left(1+\frac{2^{p+1}-1}{2^{p+1}}\right) & \text{if }s=0, e>e_{min}, k=0\\
2^e\left(1+\frac{2k-1}{2^{p+1}}\right) & \text{else}
\end{cases}
\end{align*}
and
\begin{align*}
\ceil{z}=\begin{cases}
z&\text{if }s=0, e=n, k=2^p-1\\
2^e\left(1+\frac{2k+1}{2^{p+1}}\right) & \text{else}
\end{cases}
\end{align*}
For $s=1$ note simply that
\begin{align}
\floor{-z}=-\ceil{z}\qquad\text{and}\qquad\ceil{-z}=-\floor{z}\label{eq:minz}
\end{align}
Since $u=2^{-(p+1)}$, when $s=0$ the condition $\floor{z}\leq \frac{z}{1-tu}\leq \ceil{z}$ becomes:
\begin{enumerate}
\item If $e=e_{min}, k=0$:
\begin{align*}
2^{e-1}\left(1-\frac{t}{2^{p+1}}\right)\leq 2^e\leq 2^e\left(1+\frac{1}{2^{p+1}}\right)\left(1-\frac{t}{2^{p+1}}\right)
\end{align*}
from which we get:
\begin{align}
-1\leq t\leq\frac{2^{p+1}}{2^{p+1}+1}\label{eq:trange1}
\end{align}
\item If $e>e_{min}, k=0$:
\begin{align*}
& ~2^e\left(1-\frac{1}{2^{p+2}}\right)\left(1-\frac{t}{2^{p+1}}\right)\\
\leq & ~2^e \\
\leq & ~2^e\left(1+\frac{1}{2^{p+1}}\right)\left(1-\frac{t}{2^{p+1}}\right)
\end{align*}
from which we get:
\begin{align}
-\frac{2^{p+1}}{2^{p+2}-1}\leq t\leq \frac{2^{p+1}}{2^{p+1}+1}\label{eq:trange2}
\end{align}
\item If $e=n, k=2^p-1$:
\begin{align*}
&~ 2^e\left(1+\frac{2(2^p-1)-1}{2^{p+1}}\right)\left(1-\frac{t}{2^{p+1}}\right) \\
\leq & ~ 2^e\left(1+\frac{2^p-1}{2^p}\right) \\
\leq & ~ 2^e\left(1+\frac{2^p-1}{2^p}\right)\left(1-\frac{t}{2^{p+1}}\right)
\end{align*}
from which we get
\begin{align}
-\frac{2^{p+1}}{2^{p+2}-3}\leq t\leq 0\label{eq:trange3}
\end{align}
\item Else: 
\begin{align*}
& ~2^e\left(1+\frac{2k-1}{2^p}\right)\left(1-\frac{t}{2^p}\right) \\
\leq & ~ 2^e\left(1+\frac{k}{2^p}\right)\leq \\
\leq & ~ 2^e\left(1+\frac{2k+1}{2^p}\right)\left(1-\frac{t}{2^p}\right)
\end{align*}
from which we get:
\begin{align}
-\frac{2^{p+1}}{2^{p+1}+2k-1}\leq t\leq \frac{2^{p+1}}{2^{p+1}+2k+1}\label{eq:trange4}
\end{align}
\end{enumerate}
The case for $s=1$ can then be derived from \eqref{eq:minz}:
\begin{align*}
\floor{-z}\leq \frac{-z}{1-tu}\leq \floor{-z} &~\Leftrightarrow~ -\ceil{z}\leq \frac{-z}{1-tu}\leq -\floor{z} \\
&~\Leftrightarrow\qquad \floor{z}\leq \frac{z}{1-tu}\leq \ceil{z}
\end{align*}
The bounds $t_{min}$ and $t_{max}$ of \eqref{eq:trange1}-\eqref{eq:trange4} are reached when the following relations are satisfied:
\begin{align*}
&\frac{1}{1-t_{max}u}=\frac{\ceil{z}}{z} & \frac{1}{1-t_{min}u}=\frac{\floor{z}}{z} & &\text{when } z\geq 0\\
&\frac{1}{1-t_{min}u}=\frac{\ceil{z}}{z} & \frac{1}{1-t_{max}u}=\frac{\floor{z}}{z} & &\text{when } z\leq 0
\end{align*}
It follows that 
\begin{align*}
\frac{\ceil{z}-\floor{z}}{\absv{z}}&=\left(\frac{1}{1-t_{max}u}-\frac{1}{1-t_{min}u}\right)\\
&=\frac{u(t_{max}-t_{min})}{(1-t_{max}u)(1-t_{min}u)}
\end{align*}
and thus
\begin{align}
\absv{z}u=\frac{(\ceil{z}-\floor{z})(t_{max}-t_{min})}{(1-t_{max}u)(1-t_{min}u)}\label{eq:absvzu}
\end{align}
where the coefficients $C(e,k):=\frac{(t_{max}-t_{min})}{(1-t_{max}u)(1-t_{min}u)}$ can be computed from Eqs. \eqref{eq:trange1}-\eqref{eq:trange4}.
\begin{enumerate}
\item  If $e=e_{min}, k=0$:
\[
C(e_{min},0)=\frac{2^{p+1}+1}{2^p(2^{p+1}-1)}
\]
\item If $e>e_{min}, k=0$:
\[
C(e,0)=\frac{2}{3}
\]
\item  If $e=n, k=2^p-1$:
\[
C(n,2^p-1)=\frac{3(2^{p+1}-1)}{2^{p+1}}
\]
\item Else:
\[
C(e,k)=\frac{2^p+k}{2^{p+1}}
\]
\end{enumerate}

We now use Eqs. \eqref{eq:trange1}-\eqref{eq:trange4} to express the possible values of $s,e,k$ for a given $t$. To check whether $k=0$ is possible one can simply use Eqs. \eqref{eq:trange1} and \eqref{eq:trange2} and to see if $e=n, k=2^p+1$ is possible one simply uses \eqref{eq:trange3}. For all other values of the exponent and the mantissa, \eqref{eq:trange4} gives for $t\geq 0$:
\begin{align}
2^p\left(-\frac{1}{t}-1\right)+\frac{1}{2}\leq k\leq 2^p\left(\frac{1}{t}-1\right)-\frac{1}{2}\label{eq:kfromtpos}
\end{align}
and when $t\leq 0$:
\begin{align}
2^p\left(\frac{1}{t}-1\right)-\frac{1}{2}\leq k\leq 2^p\left(-\frac{1}{t}-1\right)+\frac{1}{2}\label{eq:kfromtneg}
\end{align}
Note that when $t\in\left[-\frac{1}{2},\frac{1}{2}\right]$, all combinations of $s,e,k$ are possible, whereas for $t\in \left[-1,-\frac{1}{2}\right]\cup \left[\frac{1}{2},1\right]$ some mantissas will have to be excluded. 
We have now enough details to describe the typical distribution. For $t\in\left[-\frac{1}{2},\frac{1}{2}\right]$, \eqref{eq:errorDensity} combined with \eqref{eq:absvzu} becomes 
\footnotesize
\begin{align}
&d(t)=\sum_{z\in \F\setminus\{-\infty,0,\infty\}}f\left(\frac{z}{1-tu}\right) \frac{u\absv{z}}{(1-tu)^2}\nonumber
\\
=&\sum_{z\in \F\setminus\{-\infty,0,\infty\}}f\left(\frac{z}{1-tu}\right) \frac{1}{(1-tu)^2}\frac{(\ceil{z}-\floor{z})(t_{max}(z)-t_{min}(z))}{(1-t_{max}(z)u)(1-t_{min}(z)u)}&\nonumber
\\
=&\sum_{k=0}^{2^{p}-1}\left(\sum_{s=0}^1\sum_{e=e_{min}}^{e_{max}}f\left(\frac{z(e,s,k)}{1-tu}\right)\frac{C(e,k)(\ceil{z(e,s,k)}-\floor{z(e,s,k)})}{(1-tu)^2}\right)\nonumber 
\\
=&\sum_{k=0}^{2^p-1}\sum_{e=e_{min}+1}^{e_{max}-1}\sum_{s=0}^1 f\hspace{-2pt}\left(\hspace{-2pt}\frac{z(e,s,k)}{1-tu}\hspace{-2pt}\right)\frac{C(e,k)(\ceil{z(e,s,k)}-\floor{z(e,s,k)})}{(1-tu)^2}+ \nonumber
\\
&\hspace{-6pt} \sum_{k=0}^{2^{p}-1}\sum_{s=0}^1 f\hspace{-2pt}\left(\hspace{-2pt}\frac{z(e_{min},s,k)}{1-tu}\hspace{-2pt}\right)\frac{C(e_{min},k)(\ceil{z(e_{min},s,k)}-\floor{z(e_{min},s,k)})}{(1-tu)^2}+ \nonumber
\\
& \hspace{-6pt}\sum_{k=0}^{2^{p}-1}\sum_{s=0}^1 f\hspace{-2pt}\left(\hspace{-2pt}\frac{z(e_{max},s,k)}{1-tu}\hspace{-2pt}\right)\frac{C(e_{max},k)(\ceil{z(e_{max},s,k)}-\floor{z(e_{max},s,k)})}{(1-tu)^2}
\label{eq:pdf1/2}
\end{align}
\normalsize
Since $C(e,k)=\frac{2^p+k}{2^{p+1}}$ does not depend on $e$ when $e_{min}<e<e_{max}$ we can re-write the first term in the sum above as
\footnotesize
\[
\sum_{k=0}^{2^p-1}\frac{C(e,k)}{(1-tu)^2}\sum_{e=e_{min}+1}^{e_{max}-1}\sum_{s=0}^1 f\left(\frac{z(e,s,k)}{1-tu}\right)(\ceil{z(e,s,k)}-\floor{z(e,s,k)})
\]
\normalsize
We now make three assumptions. Each assumption is formalized by exploiting the previous one(s).

\begin{enumerate}
\item \textbf{Assumption 1.} The probability density function is constant at the scale of the intervals between floating point numbers, i.e.
\begin{align*}
&\Pro{\round(x)=z(e,s,k)}\\
=&\int_{\floor{z(e,s,k)}}^{\ceil{z(e,s,k)}} f(x)~dx \\
\approx & f\left(\frac{z(e,s,k)}{1-tu}\right)(\ceil{z(e,s,k)}-\floor{z(e,s,k)})
\end{align*}
for all values of $t$ such that $\one_{\fintvl[z]}\left(\frac{z}{1-tu}\right)=1$.
\item \textbf{Assumption 2.} The probability under $f$ of sampling a number whose rounding has exponent $e_{min}$ or $e_{max}$ is close to zero, i.e.
\footnotesize
\[
 \sum_{k=0}^{2^{p}-1}\sum_{s=0}^1 f\left(\frac{z(e_{min},s,k)}{1-tu}\right)(\ceil{z(e_{min},s,k)}-\floor{z(e_{min},s,k)})\approx 0
\]
\normalsize
and
\footnotesize
\[
 \sum_{k=0}^{2^{p}-1}\sum_{s=0}^1 f\left(\frac{z(e_{max},s,k)}{1-tu}\right)(\ceil{z(e_{max},s,k)}-\floor{z(e_{max},s,k)})\approx 0
\]
\normalsize
for all values of $t$ such that $\one_{\fintvl[z]}\left(\frac{z}{1-tu}\right)=1$
\item \textbf{Assumption 3.} When rounding a sample drawn from the distribution $f$, all mantissas are equally likely, i.e. for a given mantissa $k$
\footnotesize
\[
\sum_{e=e_{min}+1}^{e_{max}-1}\sum_{s=0}^1 f\left(\frac{z(e,s,k)}{1-tu}\right)(\ceil{z(e,s,k)}-\floor{z(e,s,k)})\approx \frac{1}{2^p}
\]
\normalsize
for all values of $t$ such that $\one_{\fintvl[z]}\left(\frac{z}{1-tu}\right)=1$.
\end{enumerate}

Under these assumptions, \eqref{eq:pdf1/2} becomes close to
\footnotesize
\begin{align}
&\sum_{k=0}^{2^p-1}\frac{C(e,k)}{(1-tu)^2}\sum_{e=e_{min}+1}^{e_{max}-1}\sum_{s=0}^1 f\left(\frac{z(e,s,k)}{1-tu}\right)(\ceil{z(e,s,k)}-\floor{z(e,s,k)})\nonumber
\\
&\approx\frac{1}{2^p(1-tu)^2}\left(\frac{2}{3}+\sum_{k=1}^{2^p-1}\frac{2^p+k}{2^{p+1}}\right)\nonumber
\\
&=\frac{1}{2^p(1-tu)^2}\left(\frac{2}{3}+\frac{3(2^{p}-1)}{4}\right)\label{eq:pdf1/2ass}
\end{align}
\normalsize
If we assume further that $(1-tu)^2\approx 1$ we get that on $\left[-\frac{1}{2},\frac{1}{2}\right]$, $d(t)$ is equal to the constant $\frac{1}{2^p}\left(\frac{2}{3}+\frac{3(2^{p}-1)}{4}\right)$ which itself is approximated well by $\frac{3}{4}$ for sufficiently large precision levels $p$.

As shown by Eqs \eqref{eq:kfromtpos} \eqref{eq:kfromtneg}, when $\absv{t}>\frac{1}{2}$, not all mantissas in the sum \eqref{eq:pdf1/2} are possible, and for $t\geq 0$ \eqref{eq:pdf1/2ass} becomes 
\begin{align*}
d(t)&\approx\frac{1}{2^p(1-tu)^2}\left(\frac{2}{3}+\sum_{k=1}^{2^p(1/t-1)-1/2}\frac{2^p+k}{2^{p+1}}\right)
\\
&=\frac{1}{2^p(1-tu)^2}\left(\frac{2}{3}+\frac{1}{2}\alpha(t)+\frac{1}{2^{p+2}}(\alpha(t)^2)\right)
\end{align*}
where $\alpha(t)=\floor{2^p(\frac{1}{t}-1)-\frac{1}{2}}$ is the usual floor function applied to $2^p(\frac{1}{t}-1)-\frac{1}{2}$. Since we cannot give an analytic expression for this integer, it is useful to consider the limit behaviour for large precisions (\ie for $p\to\infty$) in which case the expression $\alpha(t)=\floor{2^p(\frac{1}{t}-1)-\frac{1}{2}}$ is close to $2^p(\frac{1}{t}-1)-\frac{1}{2}$. The equation above then becomes
\[
d(t)\approx\frac{1}{(1-tu)^2}\left(\frac{1}{2}\left(\frac{1}{t}-1\right)+\frac{1}{4}\left(\frac{1}{t}-1\right)^2\right)
\]
in the limit where $p\to\infty$. The case where $t\leq -\frac{1}{2}$ is treated in the same way and yields the same asymptotic distribution. In the limit of large precisions and if we take $(1-tu)^2\approx 1$ we thus get the \emph{typical error distribution} with density function 
\begin{align*}
d(t)=\begin{cases}
\frac{3}{4}&\text{if }t\in\left[-\frac{1}{2} ,\frac{1}{2}\right]  \\
\frac{1}{(1-tu)^2}\left(\frac{1}{2}\left(\frac{1}{t}-1\right)+\frac{1}{4}\left(\frac{1}{t}-1\right)^2\right) & \text{else }\\
\end{cases}
\end{align*}
as given in \cref{eq:typicalpdf} and plotted in Fig. \ref{fig:typical}.

\end{document}